\documentclass[aps,pra,twocolumn,showpacs,showkeys,amsmath,amssymb,superscriptaddress]{revtex4}
\usepackage{amsthm,amsmath,amsfonts,amssymb,amsxtra,appendix,bookmark,dsfont,latexsym,stix}

\makeatletter
\usepackage{hyperref}
\usepackage{delarray,color}
\usepackage{euscript}
\usepackage{color,graphicx}

\theoremstyle{plain}

\theoremstyle{remark}

\numberwithin{equation}{section}

\DeclareMathOperator{\tr}{Tr}

\def\geqslant{\ge}
\def\leqslant{\le}
\def\bq{\begin{eqnarray}}
	\def\eq{\end{eqnarray}}
\def\bqq{\begin{eqnarray*}}
	\def\eqq{\end{eqnarray*}}

\def\eps{\varepsilon}

\def\id{\mathbb{1}}

\newcommand{\norm}[1]{\left\lVert #1 \right\rVert}

\newcommand\1{{\ensuremath {\mathds 1} }}

\newcommand{\im}{\mathrm{i}}

\renewcommand{\epsilon}{\varepsilon}
\newcommand{\bx}{\mathbf{x}}

\newcommand{\bR}{\mathbf{R}}
\newcommand{\by}{\mathbf{y}}

\def\R {\mathbb{R}}

\def\C {\mathbb{C}}
\def\N {\mathcal{N}}

\def\cE {\mathcal{E}}

\def\R {\mathbb{R}}
\def\C {\mathbb{C}}
\def\N {\mathbb{N}}

\def\gS{\mathfrak{S}}

\renewcommand{\div}{{\rm div}}
\renewcommand{\leq}{\leqslant}
\renewcommand{\geq}{\geqslant}

\newcommand{\asym}{\mathrm{asym}}
\newcommand{\curl}{\mathrm{curl}}
\newcommand{\bA}{\mathbf{A}}

\newcommand{\bp}{\mathbf{p}}
\newcommand{\bJ}{\mathbf{J}}

\newcommand{\rmd}{\mathrm{d}}

\newcommand{\cEH}{\mathcal{E}^{\rm H}}
\newcommand{\EH}{{E}^{\rm H}}
\newcommand{\gammaH}{\gamma^{\rm H}}
\newcommand{\mH}{m^{\rm H}}
\newcommand{\ELDA}{{E}^{\rm LDA}}
\newcommand{\cELDA}{\mathcal{E}^{\rm LDA}}

\newcommand{\cEcla}{\mathcal{E}^{\rm cla}}

\date{January 2026}

\begin{document} 

\title{Magnetic Thomas-Fermi theory for 2D abelian anyons}

\author{Antoine Levitt}
\affiliation{Universit\'e Paris Saclay, LMO (UMR 8628)}
\email{antoine.levitt@universite-paris-saclay.fr}

\author{Douglas Lundholm}
\affiliation{Uppsala University, Department of Mathematics}
\email{douglas.lundholm@math.uu.se}

\author{Nicolas Rougerie}
\affiliation{Ecole Normale Sup\'erieure de Lyon \& CNRS, UMPA (UMR 5669)}
\email{nicolas.rougerie@ens-lyon.fr}

\begin{abstract}
Two-dimensional abelian anyons are, in the magnetic gauge picture,
represented as fermions coupled to magnetic flux tubes. For the ground
state of such a system in a trapping potential, we theoretically and
numerically investigate a Hartree approximate model, obtained by
restricting trial states to Slater determinants and introducing 
a self-consistent magnetic field, locally proportional to matter density. 
This leads to a fermionic variant of the Chern-Simons-Schr\"odinger system. 
We find that for dense systems, a semi-classical approximation yields qualitatively good results. Namely, we derive a density functional theory of magnetic Thomas-Fermi type, which correctly captures the trends of our numerical results. In particular, we explore the subtle dependence of the ground state with respect to the fraction of magnetic flux units attached to particles.
\end{abstract}


\maketitle

 The recognition that ``anyons'' with exotic exchange statistics (different from those of bosons and fermions) were a theoretical possibility in a 2D world historically came together with the observation that they would be equivalent to the attachment of magnetic flux to ordinary particles~\cite{LeiMyr-77,Wilczek-82a,Wilczek-82b,GolMenSha-80,GolMenSha-81}. Nowadays, artificial magnetic flux attachment is the chief proposal to implement or simulate the behavior of 2D anyons, starting with the case of quasi-particles in the fractional quantum Hall effect~\cite{AroSchWil-84,LunRou-16,LamLunRou-22} and continuing with more recent proposals using artificial gauge fields in a cold atoms context~\cite{EdmEtalOhb-13,ValWesOhb-20,YakLem-18,YakEtal-19,Yakaboylu-etal-20}. Experimental evidence was found in the quantum Hall context~\cite{BarEtalFev-20,NakEtalMan-20} for 2D anyons and with cold atoms for 1D anyons~\cite{FroEtalTar-22,KwaEtalGre-24,DhaEtalNag-24} (the latter being distinct from their 2D counterparts, see~\cite{ChiEtalCel-22,RouYan-23a,RouYan-23b,RouYan-25,ValOhb-24,Yang-25,IacCabTarCel-25} for related discussions).

 An immediate difficulty is that 2D anyons, being equivalent to quite peculiar composite fermions or bosons, are very difficult to study analytically even in the simplest of situations, with no other interactions than those coming from the flux attachment (``free'' anyons). It is therefore important to be able to rely on good approximate models. Our purpose here is to propose a simple effective density functional theory for the ground state of 2D anyons in a trapping potential, inspired by the typical set-up of cold atoms experiments. This was studied before for ``almost-bosonic'' anyons~\cite{LunRou-15,CorLunRou-16,CorDubLunRou-19,Girardot-19,ValWesOhb-20,GirLee-24}, where the flux attached to bosons is assumed small in the limit of large particle numbers, and for ``almost-fermionic'' anyons~\cite{ChiSen-92,GirRou-21,GirRou-22} for a similar situation but with fermion-based anyons instead. Here we aim at a description valid for more general flux, away from the almost-bosonic and almost-fermionic limits.   
 
 We start from a fermion-based magnetic gauge picture for $N$ 2D
 anyons, i.e. from $N$ fermions coupled to magnetic flux tubes of
 intensity $2\pi \alpha \in [0,2\pi[$,  with $\alpha$ the exchange statistics parameter (counted from the fermionic end). To approximate ground-state properties at large $N$ and fixed $\alpha$, we perform a mean-field approximation. Available states are reduced to Slater determinants/quasi-free states, determined by a one-particle reduced density matrix. We solve for an approximate ground state by minimizing a Hartree functional with a self-consistent magnetic field, leading to Chern-Simons-Schr\"odinger-like variational equations. This general approach is akin to that of~\cite{HuMurRaoJai-21} (see also~\cite{HuJai-19}), where a Kohn-Sham method leads to a related model. We however fully take into account the contribution of the functional derivative of the Chern-Simons field to the variational equations. In the FQHE context, bosonic Chern-Simons-Schr\"odinger (CSS) was proposed in \cite{ZhaHanKiv-89} and fermionic CSS in \cite{HalLeeRea-93,LopFra-91}, see e.g. the discussion in~\cite[Section 5.16]{Jain-07}, as well as~\cite{TouHerHan-20} and references therein for more recent developments.

 We numerically solve the fermionic Chern-Simons-Schr\"odinger variational problem and compare the results with those derived from a semi-classical/local density approximation. The latter is an instance of a magnetic Thomas-Fermi (mTF) theory, but with a self-consistent magnetic field. 
 

 The ground-state energy and density profiles in the local density approximation depend on $\alpha$ through the expression
 \begin{align*}
c(\alpha) &:= \alpha + 2 \alpha \lfloor \alpha^{-1} \rfloor - \left(\alpha \lfloor \alpha^{-1} \rfloor \right) ^2 - \alpha ^2 \lfloor \alpha^{-1} \rfloor\\
&= 1 + \alpha^2 \left(1-\left\{ \alpha^{-1}\right\} \right)\left\{ \alpha^{-1}\right\},
 \end{align*}
where $\lfloor \alpha^{-1} \rfloor$ denotes the integer part of $\alpha^{-1}$ and 
$$ \left\{ \alpha^{-1}\right\} := \alpha^{-1} - \lfloor \alpha^{-1} \rfloor \in [0,1[$$
its fractional part;
see Figure~\ref{fig:TF constant} and compare with our numerical result in Figure~\ref{fig:num ener}.

Observe that for the more traditionally studied~\cite{Laughlin-88,FetHanLau-89,DaiLevFetHanLau-92,CheWilWitHal-89,HuMurRaoJai-21} cases of $\alpha = 1 / n$ with integer $n$ there is hence no dependence on $n$: $c(\alpha) = 1$. For more general values of $\alpha$, which do not seem to have been thoroughly considered before, the dependence on $\alpha$ of ground-state energies and density profiles is quite small. Capturing it numerically is quite challenging: we tested systems of up to $N = 100$ particles, observing a satisfying qualitative agreement with theory, improving with increasing $N$. The dependence of direct-space densities on $\alpha$ stays quite small,  leading us to suggest to look for signatures of anyon statistics in momentum-space densities instead, as put forward e.g. in~\cite{GirRou-21} in the almost-fermionic limit $\alpha \to 0$.  Such densities could be e.g. accessed by time-of-flight measurements in cold atomic gases~\cite[Section~III.A]{BloDalZwe-08}. 
 
We proceed to define the microscopic anyon and effective Hartree and CSS models more precisely in Sections~\ref{sec:anyons} and \ref{sec:Hartree}, followed by a brief recollection of the almost-fermionic limit in Section~\ref{sec:almost}. In Section~\ref{sec:mTF}, we introduce the mTF model, and in Section~\ref{sec:numerics}, we study it numerically. Finally, in Section~\ref{sec:conclusions}, we provide some concluding remarks and outlooks for future research.

\section{Magnetic gauge picture for 2D anyons}\label{sec:anyons}

We refer to the literature~\cite{Lundholm-Ency,Myrheim-99,Khare-05,Ouvry-07,Wilczek-90,LamLunRou-22} for general background, and in particular the link between (multi-valued) quantum wave functions with exchange phases $e^{\im\pi(\alpha+1)}$ and the following model. 
We consider a gas of $N$ fermions moving non-relativistically in a plane $\R^2$, each coupled to a magnetic flux tube of intensity $2\pi\alpha$, and trapped in a confining (for example harmonic) external potential 
$$ 
V:\R^2 \mapsto \R^+.
$$
Thus, the many-body Hamiltonian (defined as the Friedrichs extension of the associated quadratic form, and in units $\hbar = c = 2m = 1$)
\begin{equation}\label{eq:hamil}
H_N^{\alpha} := \sum_{j=1}^N \left[ \left( -\im \nabla_{\bx_j} + \alpha \bA (\bx_j) \right) ^2 + N V (\bx_j) \right]
\end{equation}
acts on the Hilbert space $L^2_\asym \left(\R^{2N}\right)$ of square-integrable complex-valued functions $\Psi_N(\bx_1,\ldots,\bx_N)$ that are antisymmetric w.r.t. particle exchange, 
and we have defined the magnetic vector potential felt by particle $j$
\begin{equation}\label{eq:A}
 \bA (\bx_j):= \sum_{k\neq j} \frac{(\bx_j-\bx_k)^{\perp}}{|\bx_j - \bx_k|^2}
\end{equation}
corresponding to the magnetic field 
\begin{equation}\label{eq:B}
 B (\bx_j) := \curl_{\bx_j} \bA (\bx_j) = 2\pi \sum_{k\neq j} \delta (\bx_j - \bx_k)
\end{equation}
that we treat as a scalar field since it always points perpendicular to the plane. Note the convention to mutliply $V$ by the particle number $N$ in~\eqref{eq:hamil}. Effectively, this will set the length scale of the gas as $O(1)$ for large $N$. Indeed, the first kinetic term in~\eqref{eq:hamil} typically weighs $O(N^2)$ for $N$ fermions in a bounded region (e.g. filling the Fermi disk, for $\alpha = 0$). 

Besides being the basic theoretical model for anyons in the magnetic gauge picture, one can show that the above Hamiltonian is the appropriate description for fermionic quantum impurities in a Laughlin liquid. For strong impurity/liquid interactions, the impurities couple to Laughlin quasi-holes, leading to composite objects whose dynamics is effectively governed by the above~\cite{CooSim-15,LunRou-16,LamLunRou-22,ZhaSreGemJai-14}. Another prospect is to generate the magnetic flux attachement by coupling neutral cold atoms to density-dependent artificial gauge fields~\cite{EdmEtalOhb-13,DalGerJuzOhb-11,ValWesOhb-20}.

We study the ground-state problem
\begin{multline}\label{eq:ener}
E(\alpha,N) := \inf \big\{ \left\langle \Psi_N  |H_N ^{\alpha} | \Psi_N \right\rangle_{L^2_\asym \left(\R^{2N}\right)} : \\ \Psi_N \in L^2_{\asym} \left(\R^{2N}\right), \ \norm{\Psi_N}_{L^2} = 1 \big\}, 
\end{multline}
for which we propose effective models, aiming at approximations for macroscopic systems $N\to \infty$, at fixed parameter $\alpha$ (which we may restrict to the interval $[0,1[$ by periodicity and complex conjugation). This contrasts with previous studies~\cite{ChiSen-92,GirRou-21,GirRou-22} valid in the almost-fermionic $\alpha \propto N^{-1/2}$ or almost-bosonic limits. The latter corresponds formally to $\alpha \to 1$ but is better understood~\cite{LunRou-15,Girardot-19,AtaLunNgu-24,GirLee-24,Lundholm-24} by studying the action of~\eqref{eq:hamil} on bosonic wave functions and letting $N\to \infty$, with $\alpha \propto N^{-1}.$

Although in the idealized model \eqref{eq:B} the charge-flux composites are pointlike, realistic emergent anyons do have some extent set by the underlying parameters, and we will in the following consider their local averaged/smeared distributions 
(this is well known as the ``average-field approximation'' \cite{Wilczek-90}).
Thus, we observe that, as per~\eqref{eq:B}, we are for fixed $\alpha$ and large $N$ effectively in a regime with magnetic field $B\propto N$. For usual 2D fermions in a constant magnetic field, this is precisely the regime where magnetic Thomas-Fermi theory emerges as an effective model~\cite{FouMad-19,LieSolYng-95,LieSolYng-94,LieSolYng-94b,Yngvason-91,Perice-22,PerRou-24}. This motivates us to adapt to our situation the semi-classical approximations of the aforementioned references, although we will not pursue their mathematically rigorous justification here.

\section{Hartree approximation and Chern-Simons-Schr\"odinger equations}\label{sec:Hartree}

If we ignore for a moment the many-body character of the vector potential~\eqref{eq:A} and magnetic field~\eqref{eq:B}, the energy functional to be minimized takes the form 
\begin{equation}\label{eq:mock}
\left\langle \Psi_N  |H_N ^{\alpha} | \Psi_N \right\rangle \simeq \tr \left[\left(-\im \nabla + \alpha \bA \right)^2 \gamma\right] + N \int_{\R^2} V \rho_\gamma
\end{equation}
in terms of the one-particle density matrix, obtained via a partial trace, 
$$ \gamma := N \tr_{2\to N} |\Psi_N \rangle \langle \Psi_N |,$$
and the associated one-particle density 
$$ \rho_\gamma (\bx) = \gamma (\bx;\bx),$$
where we identify $\gamma$ with its integral kernel. Minimizing the above
quadratic functional amongst fermionic wave functions leads to
a Slater determinant 
\begin{equation}\label{eq:Slater}
\Psi_N (\bx_1,\ldots,\bx_N) = \frac{1}{\sqrt{N!}}\det \left( u_i (\bx_j) \right)_{1\leq i,j \leq N},
\end{equation}
where $u_1,\ldots,u_N$ 
are orthonormal orbitals in the one-body Hilbert space $L^2(\R^2)$.
Such a quasi-free state is entirely characterized by the associated density matrix which takes the form 
\begin{align}\label{eq:gamma}
 \gamma (\bx;\by) = \sum_{j=1}^N u_j (\bx) \overline{u_j (\by)},
 \quad \text{i.e.} \quad 
 \gamma = \sum_{j=1}^N |u_j\rangle \langle u_j|.
\end{align}
We propose to approximate the original $N$-body problem by reinstating the self-consistency of the vector potential (and thus magnetic field) in~\eqref{eq:mock}, namely by setting 
\begin{equation}\label{eq:A Hartree}
 \bA (\bx) = \bA [\rho_\gamma] (\bx) := \int_{\R^2} \frac{(\bx-\by)^{\perp}}{|\bx - \by|^2}\rho_\gamma(\by) \,\rmd\by,
\end{equation}
so that the magnetic field be proportional to matter density, as in~\eqref{eq:B}:
\begin{equation}\label{eq:B Hartree}
 \curl\, \bA [\rho_\gamma] (\bx) = 2 \pi \rho_\gamma (\bx).
\end{equation}
This way we obtain the Hartree-like functional  
\begin{equation}\label{eq:Hartree}
\cEH_{\alpha} [\gamma] := \tr\left( \left(-\im \nabla + \alpha \bA [\rho_\gamma]\right)^2 \gamma \right) + N \int_{\R^2} V \rho_\gamma
\end{equation}
that we should minimize amongst all fermionic one-particle density matrices $\gamma$, leading to the energy
\begin{multline}\label{eq:Hartree inf}
\EH (N,\alpha) := \inf\big\{ \cEH_{\alpha} [\gamma] : \gamma \in \gS^1(L^2 (\R^2)), \\ 0\leq \gamma \leq \id, \ \tr \gamma = N, \ \gamma \mbox{ of rank  } N\big\},
\end{multline}
where $\gS^1(L^2 (\R^2))$ denotes the trace-class operators over $L^2
(\R^2)$, and the operator constraint $\gamma \leq \id = \mbox{Identity
  Operator}$ reflects the Pauli principle. The constraints in~\eqref{eq:Hartree inf} impose that $\gamma$ is a projector, $\gamma^2 = \gamma$.

  The first variation of the energy \eqref{eq:Hartree} with respect to an infinitesimal variation $\gamma = \gammaH + \delta \gamma$ is $\tr(H_{\gammaH} \delta\gamma)$ with 
\begin{equation}\label{eq:vareq} 
  H_{\gammaH}:= \left(-\im \nabla + \alpha \bA \left[\rho_{\gammaH}\right]\right)^2 -2 \alpha \nabla^\perp w \ast \bJ\left[\gammaH\right] + N V,
\end{equation}
where we denote 
\begin{equation*}
w(\bx) = \log |\bx|, \quad \nabla^\perp w (\bx) = \frac{\bx^\perp}{|\bx|^2},
\end{equation*}
and 
\begin{equation} 
\bJ\left[\gammaH\right] = \frac{1}{2} \sum_{j=1} ^N \left(\overline{u_j} \left(-\im \nabla + \alpha \bA \left[\rho_{\gammaH}\right] \right) u_j + \mathrm{h.c.} \right)  
\end{equation}
the current operator associated to $\gammaH$. The second term in~\eqref{eq:vareq} (where $\ast$ denotes a convolution product) is perhaps the less usual. It stems from the fact that in~\eqref{eq:Hartree} the kinetic energy operator itself depends on the density.

For a minimizer $\gammaH$ of~\eqref{eq:Hartree}, the optimality condition stemming from the first variation of the energy  leads to a self-consistent equation of the form 
\begin{equation}\label{eq:Hartree eq}
\left[\gammaH, H_{\gammaH}\right] = 0
\end{equation}
In other words, in terms of the orbitals $u_1,\ldots,u_N$ of $\gammaH$ we find 
$$H_{\gammaH} \, u_j = \lambda_j u_j$$
for eigenvalues $\lambda_j,j=1\ldots N$ of the mean-field Hamiltonian $H_{\gammaH}$. This is a generalization to several orbitals of the Chern-Simons-Schr\"odinger system. Indeed, restricting to a single orbital $N= 1$, we find an equation for a complex classical field $u\colon \R^2 \to \C$, 
$$ \left( -\im \nabla + \alpha \bA\left[|u|^2\right]\right)^2 u - 2 \alpha \nabla^\perp w \ast \bJ[u]\, u + Vu = \lambda u,$$
with the self-generated couplings
\begin{align*}
\curl \,\bA &\left[ |u|^2\right] = 2 \pi |u|^2, \\ 
\bJ[u] &= \frac{1}{2} \overline{u}\left(-\im \nabla + \alpha \bA \left[ |u|^2\right] \right)u + \mathrm{h.c.}.
\end{align*}
Moving to a time-dependent setting (i.e. replacing 
$\lambda$ with $\im \partial_t$) 
we may recognize the equations of motion for a Schr\"odinger field $u$ coupled to a Chern-Simons gauge field $\bA$ in Coulomb gauge; see e.g.~\cite{JacPi-90,Dunne-95} or~\cite[Remark~1.1]{GirLee-24}. We further mention~\cite{Ataei-24,AtaLunNgu-24,AtaGirLun-25,Lundholm-24,Nguyen-24,Ataei-25} and references therein for recent literature and connections with anyons.

\section{The almost-fermionic limit}\label{sec:almost}

For later comparison, we recall from~\cite{GirRou-21,Girardot-21} the behavior of the Hartree minimization problem~\eqref{eq:Hartree inf} in the 
almost-fermionic 
limit 
\begin{equation}\label{eq:quasiferm}
N \to \infty, \ \ \alpha = \frac{\beta}{\sqrt{N}}, \ \ \beta \mbox{ fixed}. 
\end{equation}
This limit is actually a limiting case of the magnetic Thomas-Fermi limit we discuss in Section~\ref{sec:mTF} below. We single it out because mathematically rigorous results are available in this case, providing useful comparison with the numerical results.

For $N$ 2D fermions in a region of size $O(1)$ one should keep in mind the heuristic
$\bp = -\im \nabla \propto \sqrt{N}$, so that the above scaling of parameters formally sets
all terms of~\eqref{eq:Hartree} of the same order of magnitude
$O(N^2)$. In this regime one finds that the following semi-classical version of~\eqref{eq:Hartree} on the ``position $\times$ momentum'' phase space becomes accurate:
\begin{equation}\label{eq:cla ener}
\cEcla [m] = \frac{1}{(2\pi)^2} \iint_{\R^4} \left(\left|\bp + \alpha \bA [m] (\bx) \right|^2 + N V(\bx)\right) m (\bx,\bp) \,\rmd\bx \rmd\bp,
\end{equation}
where the argument is a positive measure $m(\bx, \bp)$ on phase space, properly normalized and satisfying a form of the Pauli principle:
\begin{equation}\label{eq:contraintes Husimi}
 0 \leq m (\bx,\bp) \leq 1, \quad \iint_{\R^4} m = (2\pi)^2 N.
\end{equation}
Further, we have defined its magnetic potential
$$
\bA [m] (\bx) := \frac{1}{(2\pi)^2}\iint_{\R^4} \frac{(\bx-\by)^\perp}{|\bx-\by|^2} m (\by,\bp) \,\rmd\by \rmd\bp. 
$$
Modulo a regularization of the basic Hamiltonian \eqref{eq:hamil} (extended anyons),
this approximate model has been rigorously derived
in~\cite{GirRou-21,GirRou-22}. The results can hence serve as benchmarks
for our numerical method (see Section~\ref{sec:num quasi ferm} below).

Minimizing the functional \eqref{eq:cla ener} leads to a minimizer of the form 
\begin{equation}\label{eq:cla min}
m^{\rm cla} (\bx,\bp) = \1 \left( \left|\bp + \alpha \bA ^{\rm TF} (\bx)\right|^2 \leq 4\pi \rho^{\rm TF} (\bx) \right),
\end{equation}
where $\rho^{\rm TF}$ minimizes a bona-fide Thomas-Fermi (TF) energy:
\begin{align*} 
\cE^{\rm TF}[\rho] &:= \int_{\R^2} \left( 2\pi  \rho^2 + N V\rho \right),\\
\cE^{\rm TF}[\rho^{\rm TF}] =
E^{\rm TF} &:= \inf \left\{ \cE^{\rm TF}[\rho] : \rho \geq 0, \int_{\R^2} \rho = N \right\},
\end{align*}
and its associated magnetic potential
\begin{equation}\label{eq:TF A}
\bA^{\rm TF} (\bx) 
:= \bA[\rho^{\rm TF}]
= \int_{\R^2} \frac{(\bx-\by)^\perp}{|\bx-\by|^2} \rho^{\rm TF} (\by) \,\rmd\by.  
\end{equation}
Hence,
\begin{equation}\label{eq:TF dens}
\rho^{\rm TF} = \frac{N}{4\pi} \left( \lambda^{\rm TF} - V \right)_+
\end{equation}
for a Lagrange multiplier $\lambda^{\rm TF}$ that one determines from the mass constraint. Given $V$, one can insert~\eqref{eq:TF dens} in~\eqref{eq:TF A} to numerically compute $\bA^{\rm TF} $ and then the momentum-space density
\begin{align}\label{eq:TF mom}
t^{\rm TF} (\bp) &:= \frac{1}{(2\pi)^2} \int_{\R^2} m^{\rm cla} (\bx,\bp) \,\rmd\bx \nonumber\\
&= \frac{1}{(2\pi)^2} \int_{\R^2}  \1 \left( \left|\bp + \alpha \bA ^{\rm TF} (\bx)\right|^2 \leq 4\pi \rho^{\rm TF} (\bx) \right) \rmd\bx.
\end{align}
Note that, with a radial trap, the TF density is also radial and Newton's theorem gives 
\begin{equation}\label{eq:TF A bis}
\bA^{\rm TF} (\bx) = \frac{2\pi\bx^\perp}{|\bx|^2} \int_{r= 0} ^{|\bx|} \rho^{\rm TF} (r) r \rmd r. 
\end{equation}
For the particular choice $V(\bx) = |\bx|^s, s>0$ one finds 
$$ \lambda^{\rm TF} = \left(\frac{4(s+2)}{s}\right)^{\frac{s}{s+2}},$$
and thus 
\begin{equation}\label{eq:TF A ter}
\bA^{\rm TF} (\bx) = \begin{cases}
                      \frac{N\bx^\perp}{4} \left(\lambda^{\rm TF} - \frac{2}{s+2} |\bx|^{s} \right) \mbox{ for } |\bx|^s \leq \lambda^{\rm TF},\\
                      \frac{N\bx^\perp}{|\bx|^2}\mbox{ for } |\bx|^s \geq \lambda^{\rm TF}.
                     \end{cases}
\end{equation}
As is usual with ordinary Thomas-Fermi theory (see also \cite{ChiSen-92}), the
position-space density~\eqref{eq:TF dens} is independent of the magnetic field, and thus of $\alpha$. Instead, the momentum-space density~\eqref{eq:TF mom} should be used to look for signatures of anyon statistics, i.e. of the Chern-Simons-like gauge field. Indeed,
\begin{equation}\label{eq:TF mom scaled}
\tilde{t} ^{\rm TF}(\bp):= N^{-1} t^{\rm TF} \left(\sqrt{N}\bp\right) 
\end{equation}
converges in the regime~\eqref{eq:quasiferm} to a non-trivial profile depending on the parameter $\beta$.

\section{Magnetic Thomas-Fermi approximation}\label{sec:mTF}

We turn to an approximation that should be valid more generally for the trapped anyon gas, namely the high-density limit 
\begin{equation}
 N \to \infty, \ \alpha \mbox{ fixed,}
\end{equation}
corresponding to $\beta \sim \sqrt{N} \to \infty$ in \eqref{eq:quasiferm}. This will be of magnetic Thomas-Fermi type, and just as usual fermionic Thomas-Fermi can be recovered as a limit case of mTF~\cite{LieSolYng-95,LieSolYng-94,LieSolYng-94b}, the almost-fermionic description of the preceding section is a limit case of the considerations below. Indeed, one can go from mTF theory to the almost-fermionic limit smoothly by letting $\alpha \to 0$ in the former (observe that from~\eqref{eq:TF constant} we clearly have that $c(\alpha) \underset{\alpha \to 0}{\to} 1$). 

We first describe the limit model briefly before presenting two parallel derivations in Sections~\ref{sec:LDA} and~\ref{sec:measures} below. After eliminating momentum degrees of freedom similarly as
above, the effective energy functional becomes 
\begin{equation}\label{eq:TF sem} 
\cE^{\rm mTF}_{\alpha} [\rho] := \int_{\R^2} \left(2 \pi c (\alpha) \rho(\bx)^2 + N V (\bx) \rho(\bx)\right)\rmd\bx,
\end{equation}
to be minimized under the constraints
$$ 
\rho \geq 0, \quad \int_{\R^2} \rho = N.
$$
Here 
\begin{align}\label{eq:TF constant}
 c(\alpha) &= \alpha + 2 \alpha \lfloor \alpha^{-1} \rfloor - \left(\alpha \lfloor \alpha^{-1} \rfloor \right) ^2 - \alpha ^2 \lfloor \alpha^{-1} \rfloor\nonumber \\
 &= 1 + \alpha^2 \left(1-\left\{ \alpha^{-1}\right\} \right)\left\{ \alpha^{-1}\right\} 
\end{align}
can be interpreted as $(2\pi)^{-1}\rho^{-2}$ times the energy per particle of the homogeneous anyon gas; see~\eqref{eq:app loc ener} below. 

We plot $c(\alpha)-1$ in Figure~\ref{fig:TF constant}, highlighting the quite non-trivial dependence on $\alpha$ and the reduction to the value for free fermions $c(\alpha) = 1$ for $\alpha=n^{-1}$, $n$ integer or when $\alpha \to 0$ (this is the almost-fermionic limit of Section~\ref{sec:almost}).

 \begin{figure}
\begin{center}
 \includegraphics[width=7cm]{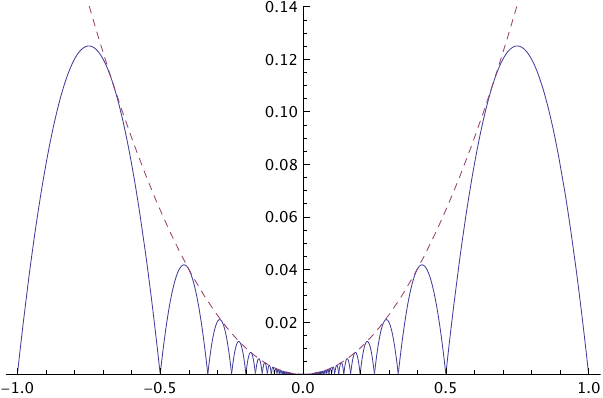}
\caption{Plot $\alpha \mapsto c(\alpha) - 1$ of the anyonic magnetic Thomas-Fermi constant \eqref{eq:TF constant}, relative to the free fermion case, i.e. fermions at $\alpha=0$ and bosons at $\alpha=1$. Also shown (dashed) is the upper bound $\alpha^2/4$.}
\label{fig:TF constant}
\end{center}
 \end{figure}

If, for concreteness, the external potential is set to be 
\begin{equation}\label{eq:pot}
 V(\bx) = |\bx|^s, \quad s >0,
\end{equation}
then explicit minimization of~\eqref{eq:TF sem} as above leads to the $\alpha$-dependent densities and energies
\begin{align}
 \rho^{\rm mTF}_{\alpha} (\bx) &= \frac{N}{4\pi c(\alpha)} \left( \lambda_{\alpha} - |\bx| ^s\right)_+, \nonumber\\
 \lambda_{\alpha} &= \left( 4c(\alpha) \frac{(s+2)}{s}\right)^{\frac{s}{s+2}}, \nonumber\\
 E^{\rm mTF}_{\alpha}&= N^2 \frac{s}{8c(\alpha)(s+1)} \lambda_{\alpha}^{\frac{2s + 2}{s}} \nonumber\\
 &=N^2 \frac{s}{8(s+1)} \left(4 \frac{s+2}{s} \right) ^{\frac{2s+2}{s+2}} c(\alpha) ^{\frac{s}{s+2}}.\label{eq:TF min}
\end{align}
We compare these expressions to numerical simulations in Section~\ref{sec:num mTF}. We first derive them using the local density approximation, and then return to the question of the distribution in momentum space using magnetic semi-classical measures. Thus, in brief, the present section aims at supporting the magnetic Thomas-Fermi approximation in the form of the following claim:
\begin{equation}\label{eq:result ener sem}
\boxed{\mbox{For fixed $\alpha$, } \EH (N,\alpha) \underset{N\to \infty}{\simeq}  E^{\rm mTF}_{\alpha}.}
\end{equation}
Moreover, with $\gammaH$ a minimizer for~\eqref{eq:Hartree inf}, the spatial density 
$$\rho^{\rm H}(\bx):=\gammaH(\bx,\bx)$$
should satisfy 
\begin{equation}
\label{eq:result dens sem}
\boxed{\rho^{\rm H} \underset{N\to \infty}{\simeq} \rho^{\rm mTF}_{\alpha}.}
\end{equation}

\subsection{Local density approximation}\label{sec:LDA}

The most direct way to argue for~\eqref{eq:result ener sem} is via a local density approximation (LDA), which we expect to be valid because the magnetic field in the Hartree model is a local function of the density. This was previously discussed briefly in~\cite[Section~III.D]{Lundholm-Ency}.

Consider a bounded reference domain $\Lambda\subset{\R^2}$ of unit area, and the thermodynamic limit of the energy density for the corresponding homogeneous problem at fixed average density $\varrho = N/L^2$:
\begin{equation}\label{eq:loc ener}
 e (\alpha,\varrho) = \lim_{N,L\to \infty,NL^{-2}= \varrho} \frac{1}{L^2} \EH (N,\alpha,L\Lambda),
\end{equation}
where for any subset $\Omega\subset \R^2$ we define the local Hartree energy analogously to \eqref{eq:Hartree}-\eqref{eq:Hartree inf}:
\begin{multline*}
\EH (N,\alpha,\Omega) := \inf\big\{ \cEH_{\alpha,\Omega} [\gamma] : \gamma \in \gS^1(L^2 (\Omega)), \\ 0\leq \gamma \leq 1, \ \tr \gamma = N, \ \gamma \mbox{ of rank } N\big\},
\end{multline*}
with
$$
\cEH_{\alpha,\Omega} [\gamma] := \tr_{L^2 (\Omega)}\left( \left(-\im \nabla + \alpha \bA [\rho_\gamma]\right)^2 \gamma \right). 
$$
Clearly, by scaling extensive quantities, we will have 
$$ e(\alpha,\varrho) \propto \varrho^2,$$
whence the point of scaling $V$ in~\eqref{eq:cELDA} below, to impose a strong density limit.

We shall argue that, with $c(\alpha)$ as in~\eqref{eq:TF constant}, we have
\begin{equation}\label{eq:app loc ener}
 \boxed{e (\alpha,\varrho) \simeq 2\pi  c(\alpha) \varrho^2.}
\end{equation}
Our claim~\eqref{eq:result ener sem} can now be rephrased as
\begin{equation}\label{eq:LDA}
 \EH (N,\alpha) \simeq \ELDA (N,\alpha),
\end{equation}
with 
\begin{equation}\label{eq:ELDA}
\ELDA (N,\alpha) := \inf\left\{ \cELDA[\rho] : \rho \geq 0, \ \int_{\R^2} \rho = N \right\} 
\end{equation}
and
\begin{equation}\label{eq:cELDA}
\cELDA[\rho] := \int_{\R^2} \bigl( e(\alpha,\rho(\bx)) + N V (\bx) \rho(\bx) \bigr) \,\rmd\bx. 
\end{equation}
Our discussion in favor of~\eqref{eq:result ener sem} therefore now relies entirely on arguing for~\eqref{eq:app loc ener}, which we now proceed to do.

For simplicity, we discuss a fixed box of side length $= 1$ with $N=\varrho\gg 1$ particles inside. We can reduce to this case by scaling length units. We assume that the density will be essentially homogeneous:
$$ 
\rho (\bx) \simeq \varrho \equiv \mathrm{constant},
$$
so that the magnetic field becomes
\begin{equation}\label{eq:B field}
B = \alpha \,\curl \, \bA [\rho] = 2\pi \alpha \varrho \equiv \mathrm{constant} \gg 1.
\end{equation}
Thus, we think of the problem at hand as equivalent to that of minimizing the energy of $\varrho$ free fermions in Landau's Hamiltonian
\begin{equation}\label{eq:Landau} 
H_B = \left(-\im \nabla + \bA \right) ^2 \mbox{ with } \curl \, \bA = B \equiv \mathrm{constant}.
\end{equation}
The energy levels of the above are essentially the Landau levels 
\begin{equation}\label{eq:LL}
E_k = 2B \left(k + \frac{1}{2}\right), \ 
k = 0,1,2,\ldots
\end{equation}
(recall that $\hbar = e = c = 1$ but the mass $m = 1/2$). These have degeneracy proportional to area times 
\begin{equation}\label{eq:LL dege}
\frac{B}{2\pi} = \alpha \varrho. 
\end{equation}
Indeed, the magnetic length of the field is, in our units,
\begin{equation}\label{eq:mag length}
\ell_B := \frac{1}{\sqrt{B}} \ll 1 = \mbox{ size of the box,}
\end{equation}
which will ensure that we may safely approximate the energy levels by those of the full plane model, modulo taking the degeneracy~\eqref{eq:LL dege} into account. Reviews of these facts are found e.g. in~\cite[Sections~2 and~3]{Perice-22},~\cite[Section~3]{Jain-07} or~\cite{RouYng-19}.

We set $n:=\lfloor \alpha^{-1}\rfloor$ and write
$$
\varrho = n \alpha \varrho + \varrho^{\rm top}, 
\qquad \varrho^{\rm top} := \left(1 - n \alpha\right)\varrho = \left\{\alpha^{-1}\right\}\alpha\varrho,
$$
to see that, taking the degeneracy~\eqref{eq:LL dege} into account, we minimize the energy by completely filling the $n$ lowest Landau levels (with $\alpha \varrho$ particles per level) and putting the remaining 
$\varrho^{\rm top} < \alpha \varrho $ 
fermions in a last, partially filled, level. 
Hence, an integer number $n$ of levels are precisely filled if $\varrho^{\rm top}=0$, i.e. if and only if the fractional part $\left\{\alpha^{-1}\right\}=0$.

Now, summing the level energies~\eqref{eq:LL} times the number of fermions occupying them, we find 
\begin{align*} 
e(\alpha,\varrho) &= \sum_{k=0}^{n-1} 2B \left(k + \frac{1}{2}\right) \alpha \varrho + 2B \left( n + \frac{1}{2}\right)\left(1 - n \alpha\right)\varrho\\
&= \sum_{k=0}^{n-1} 4\pi \alpha^2 \varrho^2 \left(k + \frac{1}{2}\right)  + 4\pi \alpha \varrho^2 \left( n + \frac{1}{2}\right)\left(1 - n \alpha\right)\\
&= 2\pi \alpha \varrho^2+ 4\pi \alpha^2 \varrho^2 \sum_{k=0}^{n-1}k + 4\pi \alpha \varrho^2 \left( 1-n\alpha\right)n\\
&= 2\pi \varrho^2 \left( \alpha + n (n-1)\alpha^2 + 2\alpha n - 2\alpha^2 n ^2 \right),
\end{align*}
and thus the claim~\eqref{eq:app loc ener} follows.

\subsection{Semi-classical measures}\label{sec:measures} 

We turn to a more semi-classical point of view on our Claim~\eqref{eq:result ener sem}. This will explain in which sense magnetic Thomas-Fermi theory should be expected to predict the behavior of  Hartree minimizers. Indeed, implicit in the local density approximation above, is a magnetic semi-classical phase space~\cite{FouMad-19,LieSolYng-95,LieSolYng-94,LieSolYng-94b,Yngvason-91,Perice-22,PerRou-24} of the form ``position of cyclotron orbit $\times$ Landau level index'', instead of the more standard ``position $\times$ momentum'' of Section~\ref{sec:almost} above. We first recall this in the context of free fermions in a large, constant, magnetic field, before turning to adaptations to the present context of 2D anyons. The general picture is derived from the well-known~\cite{ChaFlo-07,RouYng-19,Goerbig-09} classical motion in a large magnetic field, split into a cyclotron motion (whose radius, quantum mechanically, is quantized proportionally to $\sqrt{n}$, with $n$ integer) and a drift of the orbit center $\bR$. 

For free fermions in a constant magnetic field, the one-body Hamiltonian \eqref{eq:Landau} in the symmetric gauge is
\begin{equation}\label{eq:Landau op}
H_B = \left( -\im \nabla_{\bx} + \bA (\bx) \right) ^2, \quad \bA (\bx) = \frac{B}{2} \bx^\perp.  
\end{equation}
One can define a coherent operator by its integral kernel 
\begin{equation}\label{eq:coh op}
\Pi_{n,\bR} ^B (\bx,\by) := g_{\eps} (\bx - \bR) \, \Pi_n ^B (\bx,\by) \, g_{\eps} (\by - \bR),
\end{equation}
where 
\begin{align}\label{eq:nLL}
\Pi_n ^B (\bx,\by) &= \frac{B}{2\pi} \exp\left(\im B (x_1 y_2- x_2 y_1) - B\frac{|\bx - \by|^2}{4}\right)\nonumber\\
&\times L_n \left(B\frac{|\bx-\by|^2}{2}\right), 
\end{align}
in terms of Laguerre polynomials
$$
L_n \left(B\frac{|\bx-\by|^2}{2}\right) = \sum_{k= 0}^n {n \choose k} \frac{(-1)^k}{k!} \left(B\frac{|\bx-\by|^2}{2}\right)^k,
$$
is the integral kernel of the orthogonal projector on the $n$-th Landau level (with energy~\eqref{eq:LL}), and $g_\eps$ is a localization function normalized in $L^2$, e.g.
\begin{equation}\label{eq:gaussian}
 g_{\eps} (\bx) = c_\eps e^{-|\bx| ^2 / \eps^2}. 
\end{equation}
One typically asks
\begin{equation}\label{eq:lengths}
B^{-1/2} \ll \eps \ll 1,
\end{equation}
i.e. $\eps$ is larger than the magnetic length but smaller than the typical system size. One should think of the operator $\Pi_{n,\bR}$ as localizing close to the $n$-th Landau level and close to the position $\bR$ (i.e., classically, close to a cyclotron orbit of radius $\propto\sqrt{n}$ around the guiding center $\bR$). 

Alternatively, one can define 
\begin{equation}\label{eq:coh op 2}
\Pi_{n,\bR} ^B (\bx,\by) := \frac{B}{2\pi}\psi_{n,\bR} ^B (\bx) \overline{\psi_{n,\bR} ^B (\by)},
\end{equation}
with the vortex coherent state ($z$ and $Z$ are complex numbers corresponding to $\bx$ and $\bR$ respectively)
\begin{equation}\label{eq:coh stat}
\psi_{n,\bR} ^B (\bx) := \frac{\im ^n \sqrt{B}}{\sqrt{2\pi n!}} \left(  \sqrt\frac{{B}}{{2}} (z-Z) \right) ^n e^{-B\frac{|\bx-\bR|^2}{4} -\frac{B\im}{2}\bR^\perp \cdot \bx},
\end{equation}
but we will mostly rely on~\eqref{eq:coh op} in the sequel.

Given a one-particle density matrix for $N$ fermions
\begin{equation}\label{eq:Pauli q}
 0 \leq \gamma \leq \1, \ \tr \gamma = N,
\end{equation}
we can form the associated Husimi function
\begin{equation}\label{eq:Husimi}
 m(n,\bR):= \tr \left[\Pi_{n,\bR}^B \, \gamma \right],
\end{equation}
with $\Pi_{n,\bR}^B$ set as either~\eqref{eq:coh op} or~\eqref{eq:coh op 2}.  Since
$$\Pi_n^B (\bx,\bx) \equiv \frac{B}{2\pi},$$
one deduces from~\eqref{eq:Pauli q} a semi-classical version of the Pauli principle:
\begin{equation}\label{eq:Pauli c}
 0 \leq m (n,\bR) \leq \frac{B}{2\pi},
\end{equation}
and since 
$$ \sum_{n \ge 0} \Pi_n^B (\bx,\by) = \delta(\bx-\by),$$
it follows that
\begin{equation}\label{eq:Pauli c 2}
\sum_{n \ge 0} \int_{\R^2} m(n,\bR) \,\rmd\bR = N. 
\end{equation}
In view of the localization properties of $\Pi_{n,\bR} ^B$ on the $(n,\bR)\in \N \times \R^2 $ phase space, one approximates the quantum energy by the semiclassical one in the manner 
\begin{multline}\label{eq:semiclass}
\cE^{q} [\gamma] := \tr\left( \left(H_B + N V\right) \gamma \right) \\ 
\simeq \sum_{n \ge 0} \int_{\R^2} \left( 2B \left(n+\frac{1}{2}\right) + N V(\bR)\right) m(n,\bR) \,\rmd \bR =: \cE^c [m].
\end{multline}
In minimizing the right-hand side of the above, one can eliminate the Landau level variable $n$ by setting
\begin{equation}\label{eq:dens class}
\varrho (\bR) := \sum_{n \ge 0} m(n,\bR) \, \simeq \, \gamma (\bR,\bR). 
\end{equation}
%
Solving next for $\varrho$ as previously, and defining
$$ n_\star (\bR) := \left \lfloor\frac{2\pi \varrho(\bR)}{B}\right\rfloor$$
the number of filled levels, one finds
\begin{equation}\label{eq:local LL dist}
m (n,\bR) = \begin{cases}
	\frac{B}{2\pi}, & \mbox{ for } 0 \leq n \leq n_\star (\bR) - 1,\\
	\frac{B}{2\pi} \left( \frac{2\pi \varrho(\bR)}{B} - n_\star (\bR)\right), & \mbox{ for } n = n_\star(\bR),\\
	0, & \mbox{ for } n > n_\star (\bR).
	\end{cases} 
\end{equation}

In the context of 2D anyons in the Hartree approximation of Section~\ref{sec:Hartree}, we follow the previous steps mutatis mutandis. We only replace any occurrence of the fixed, constant, magnetic field $B$ by 
\begin{equation}\label{eq:BR}
 B (\bR) := 2\pi \alpha \rho_\gamma (\bR),
\end{equation}
with $\rho_\gamma$ the system's spatial density, in accordance with~\eqref{eq:B Hartree}.

We start from the Hartree energy
$$
\cEH_{\alpha} [\gamma] := \tr\left( \left(\widetilde{H}_\gamma + N V \right) \gamma \right) 
$$ 
with 
$$\widetilde{H}_\gamma := \left(-\im \nabla + \alpha \bA [\rho_\gamma]\right)^2$$
and a minimizer $\gammaH$ thereof, with density $\rho^{\rm H}$. We form the Husimi function in the manner 
\begin{equation}\label{eq:Husimi any}
\mH (n,\bR) := \tr\left[ \Pi_{n,\bR} ^{B(\bR)} \, \gammaH \right],
\end{equation}
with $\Pi_{n,\bR} ^{B(\bR)}$ as in~\eqref{eq:coh op} or~\eqref{eq:coh op 2}, and 
$$B(\bR)\simeq 2\pi \alpha \sum_{n \ge 0} \mH (n,\bR),$$ in accordance with~\eqref{eq:BR}. 
We then expect $\cEH_{\alpha} [\gamma]$ to be approximated by
\begin{align}\label{eq:semiclass any}
\cE^{c}_\alpha [m] &:= \sum_{n \ge 0} \int_{\R^2} \left( 4\pi \alpha \rho_\gamma (\bR) \left(n+\frac{1}{2}\right) + N V(\bR)\right) m(n,\bR) \,\rmd \bR, \nonumber\\
\rho_\gamma (\bR) &= \sum_{n \ge 0} m(n,R),
\end{align}
which should be minimized under the constraints 
$$ m (n,\bR) \leq \alpha \rho_\gamma (\bR), \qquad 
\sum_{n \ge 0} \int_{\R^2} m(n,R) \,\rmd \bR = N.$$
An explicit solution leads to a literal translation of~\eqref{eq:local LL dist}:
\begin{equation}\label{eq:local LL dist any}
\mH (n,\bR) \simeq \rho^{\rm H} (\bR) \times \begin{cases}
	\alpha, & \mbox{ for } 0 \leq n \leq \left\lfloor \alpha^{-1} \right\rfloor - 1,\\
	\alpha \left\{ \alpha^{-1} \right\},
	& \mbox{ for } n = \left\lfloor \alpha^{-1} \right\rfloor, \\
	0, & \mbox{ for } n > \left\lfloor \alpha^{-1} \right\rfloor.
	\end{cases} 
\end{equation} 
The dependence on $n$ in the above should be thought of as the distribution within local Landau levels leading to the energy~\eqref{eq:TF sem}. Approximating $\rho^{\rm H}$ by $\rho^{\rm TF}_{\alpha}$ gives an explicit formula, and the final proposed approximation for the full ground-state density matrix is 
\begin{equation}\label{eq:up symb}
\gammaH \simeq \sum_{n\geq 0} \int_{\R^2} \mH (n,\bR) \,\Pi_{n,\bR} ^{B(\bR)} \,\rmd \bR,
\end{equation}
from which one can extract relevant physical quantities.

\section{Numerical solutions}\label{sec:numerics}

We turn to the direct numerical minimization of the Hartree energy functional \eqref{eq:Hartree}-\eqref{eq:Hartree inf} from Section~\ref{sec:Hartree}, and the comparison of the results with the expectations derived in Sections~\ref{sec:almost} and~\ref{sec:mTF}. An approach akin to ours has been compared to other approximation schemes (exact results and trial wave-functions approach) in~\cite{HuMurRaoJai-21}. Comparison to e.g. exact diagonalization or Density Matrix Renormalization Group approaches would be useful as well, see e.g. ~\cite{HuJai-19} and~\cite[Section 5.16]{Jain-07} for such comparisons in the FQHE context at small particle numbers. They seem difficult in our context, for we need to consider quite large particle numbers (typically $N\sim 100$, which demands a very refined numerical resolution, see Section~\ref{sec:nummeth} below) to observe clear signatures. Indeed the $\alpha$-dependent features that we observe are suggested by an approximate theory which is expected to become exact in the $N\to \infty$ limit. Our best benchmark is therefore the comparison of the numerically obtained results with the exact known ones~\cite{GirRou-21,GirRou-22} in the almost-fermionic limit. The latter (for small $\alpha$) are much more theoretically grounded than those for general $\alpha$ involving mTF. We hence think of their reproduction by the simulations in Section~\ref{sec:num quasi ferm} below as a testbed for the numerical method, giving confidence for the more challenging simulations of Section~\ref{sec:num mTF}.

\subsection{Methods}\label{sec:nummeth}
We minimize the functional \eqref{eq:Hartree} as a
function of the orbitals $u_{j}$, $j= 1,2,\ldots,N$ using a small adaptation of the code
DFTK \cite{HerLevCan-21}. This code was originally designed to solve the Kohn-Sham equations for
electrons, but has an architecture flexible enough to accommodate
\eqref{eq:Hartree} easily. We discretize the orbitals using a spectral
method in a large box with periodic boundary conditions. We use a standard
Riemannian L-BFGS method on the Stiefel manifold to perform the
minimization \cite{Boumal-23}.

The numerically most sensitive part of the scheme is the solution of
the system
$$\curl \, \bA = 2\pi \rho, \quad \div \, \bA = 0, \quad \rho= \sum_{j=1}^N |u_j|^2.$$
Solving this equation naively in the Fourier domain imposes an artificial periodicity of $\bA$ which is incompatible with its slow decay at infinity, leading to very slow
convergence with respect to the box size. This problem is strictly
analogous to that encountered in density functional theory simulations of electrons in
charged cells (e.g. ions), where the Coulomb potential is long-ranged.
To correct for this, we analytically solve the equation
$$\curl\, \bA_{\rm ref} = 2\pi \rho_{\rm ref}$$
for a gaussian $\rho_{\rm ref}$ of reasonable width and of total integral equal to
that of $\rho$. Then, we solve for the correction
$$\curl \, \bA_{\rm diff} = 2\pi (\rho - \rho_{\rm ref})$$ by a Fourier
method, assuming periodic boundary conditions for $\bA_{\rm diff}$, and
reconstruct the total potential 
$$\bA = \bA_{\rm ref} + \bA_{\rm diff}$$ in real
space. Since $\rho$ is expected to be roughly spherical,
$A_{\rm diff}$ decays quickly at infinity, and this method is found to
yield rapid convergence with respect to the box size in practice.

The Thomas-Fermi solution \eqref{eq:TF dens} can be used to provide
rough guesses of the box and grid size. In particular, in the scaling
we use, the width of the Thomas-Fermi solution does not scale with
$N$, but its momentum-space width does. We use this to guide the
selection of numerical parameters in our numerical experiments, by
taking a box size equal to $2 q_{x}$ times the Thomas-Fermi real-space
width, and a grid size equal to $2\pi q_{p}$ over the momentum-space
width (which scales with $N$). The parameters $q_{x}$ and $q_{p}$ are
quality factors. Due to the spectral accuracy of our numerical scheme
and the strong locality of $\rho$ and $A_{\rm diff}$, the convergence
with respect to both $q_{x}$ and $q_{p}$ is very fast; we found
$q_{x} = 2, q_{p} = 6$ to give visually converged results for the
values of $N$ reported. 

A further key technical point concerns the accuracy of the radially
averaged densities we plot below. It is guaranteed by performing the
averaging exactly, i.e. using the analytical expression for the
average of a plane wave (the average over a circle of radius $r$ of
$e^{iq\cdot x}$ is given in terms of Bessel functions as $J_{0}(|q|\,r)$).

\subsection{Almost-fermionic regime}\label{sec:num quasi ferm}

We first check the agreement between our numerical method and the theoretical findings recalled in Section~\ref{sec:almost}. Some preliminary checks appeared in~\cite[Section~1.5.4]{Girardot-21}. Our new simulations confirm those.

We solve for $N$ orthonormal orbitals $u_1, \ldots, u_N$ and form the associated density matrix
$$ \gamma (\bx,\by) = \sum_{j=1}^N u_j (\bx) \overline{u_j(\by)}.$$
In the almost-fermionic limit, the densities in position space 
\begin{equation}\label{eq:dens pos}
\rho(\bx) = \gamma (\bx,\bx) = \sum_{j=1} ^N |u_j (\bx)| ^2 
\end{equation}
and momentum space (with $\widehat{u}$ the Fourier transform of $u$)
\begin{equation}\label{eq:dens mom}
t(\bp) = \sum_{j=1} ^N |\widehat{u_j} (\bp)| ^2 
\end{equation}
should both be correctly approximated by the Thomas-Fermi description
from Section~\ref{sec:almost}. We plot their average over rotations after
dividing by $N$, and rescaling $\bp$ by a factor $N^{-1/2}$ in momentum
space (recall that~\eqref{eq:TF mom scaled} is expected to
converge for large $N$). Comparisons are then made with the explicit~\eqref{eq:TF
  dens} (insert $c(\alpha)= 1$ in~\eqref{eq:TF min}) and
semi-explicit~\eqref{eq:TF mom} (where the integrations involving the
TF position density are performed numerically using a Monte Carlo
method \cite{Lepage-78,Hahn-05}). Recall that the expected position density is independent of $\alpha$ in this regime. For reference we also plot the momentum density one would obtain from ordinary TF theory, i.e. at $\alpha = 0$.

We provide two typical sets of plots:

 \begin{figure}
\begin{center}
 \includegraphics[width=8cm]{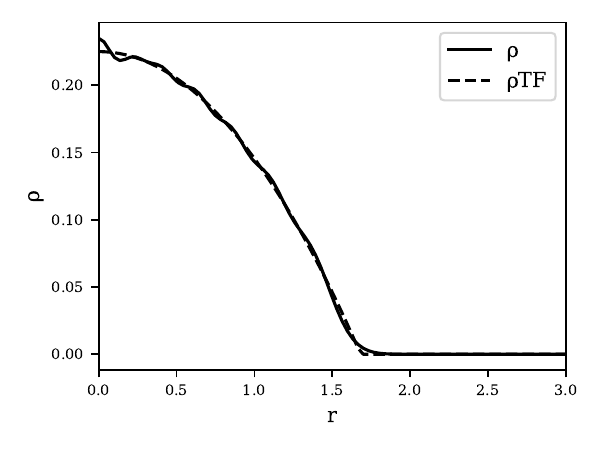}
 \caption{Position density (normalized), $N=100$, $\alpha=\frac{1.5}{\sqrt{N}}=0.15$,
   $V(\bx)=|\bx|^{2}$. }
\label{fig:almost 1}
\end{center}
 \end{figure}

 \begin{figure}
\begin{center}
 \includegraphics[width=8cm]{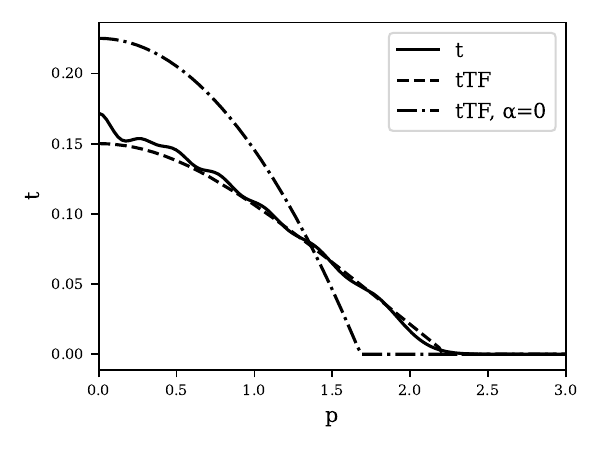}
 \caption{Momentum density (rescaled), $N=100$, $\alpha=\frac{1.5}{\sqrt{N}}=0.15$,
   $V(\bx)=|\bx|^{2}$. }
\label{fig:almost 1 bis}
\end{center}
 \end{figure}

\begin{figure}
\begin{center}
 \includegraphics[width=8cm]{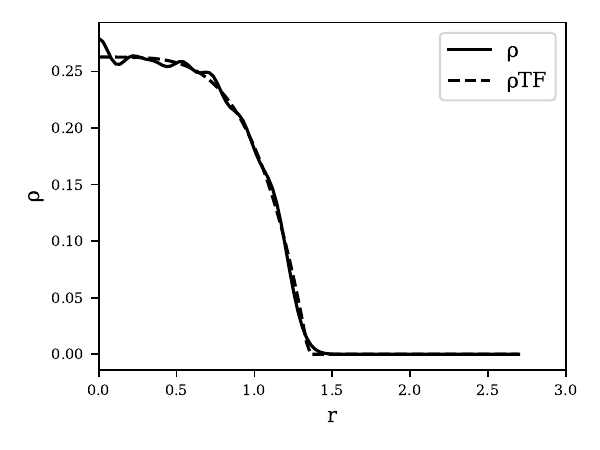}
 \caption{Position density (normalized), $N=100$, $\alpha=\frac{1.5}{\sqrt{N}}=0.15$,
   $V(\bx)=|\bx|^{4}$. }
\label{fig:almost 2}
\end{center}
 \end{figure}

\begin{figure}
\begin{center}
 \includegraphics[width=8cm]{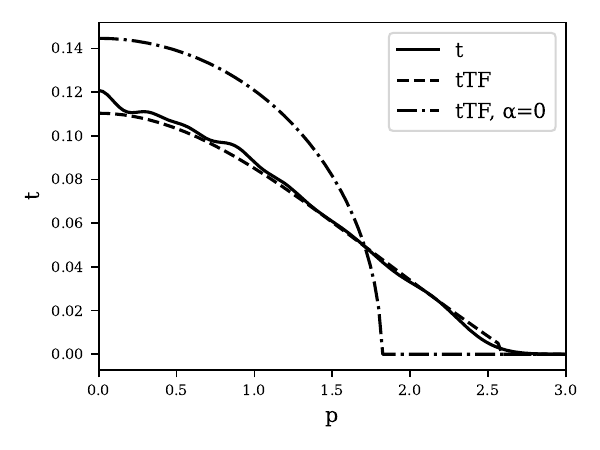}
 \caption{Momentum density (rescaled), $N=100$, $\alpha=\frac{1.5}{\sqrt{N}}=0.15$,
   $V(\bx)=|\bx|^{4}$. }
\label{fig:almost 2 bis}
\end{center}
 \end{figure}

\begin{itemize}
 \item Figures~\ref{fig:almost 1} and~\ref{fig:almost 1 bis} correspond to $N= 100$ particles in a harmonic trap 
\begin{equation}\label{eq:harm pot}
 V(\bx) = |\bx|^2 
\end{equation}
with 
$$ \alpha = 1.5 N^{-1/2} = 0.15.$$
Note that in a harmonic trap, the usual TF momentum and position densities (i.e. at $\alpha = 0$) are identical.  
\item Figures~\ref{fig:almost 2} and~\ref{fig:almost 2 bis} correspond to $N= 100$ particles in a quartic trap 
$$ V(\bx) = |\bx|^4$$
with again $\alpha = 0.15$.
\end{itemize}

In both cases we find a rather good agreement between the theoretical expectations and the numerical simulations, in particular in view of the fact that the numerical densities result from averaging a sum of many rather oscillating profiles over rotations. The agreement found already occurs for a ``moderately small'' value of $\alpha$, indicating a rapid convergence of the Hartree minimization problem to the almost-fermionic Thomas-Fermi one for small $\alpha$. This could be guessed from the rapid convergence 
$$c(\alpha)\underset{\alpha \to 0}{\to} 1$$
apparent in Figure~\ref{fig:TF constant}, and is in agreement with previous simulations~\cite{HuMurRaoJai-21}.  As anticipated, the position-space density is indistinguishable from the purely fermionic case $\alpha = 0$, $c(\alpha) = 1$. The momentum-space density is much more interesting, and markedly differs from that expected for free fermions, thus showing a clear signature of anyonic statistics via the interaction with the Chern-Simons-like gauge field.

\subsection{Magnetic Thomas-Fermi regime}\label{sec:num mTF}

We now move away from the almost-fermionic limit to check the prediction of mTF theory at general $\alpha$. Inserting reasonable numbers in~\eqref{eq:TF min} immediately reveals that, as regards energies and position-space densities, we are looking for relatively small effects. Roughly (cf. again the relatively small variations in Figure~\ref{fig:TF constant}) at the level of energies, one looks for variations in the $5\%$ range. We thus have to perform rather challenging simulations (in terms of particle number, grid size and computation time). We estimate the precision of our numerical results as circa $1\%$.   

\begin{figure}
\begin{center}
 \hspace{-1cm}\includegraphics[width=8cm]{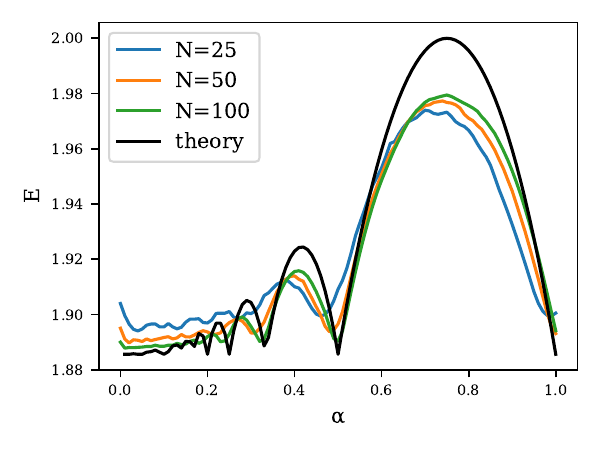}
\caption{Ground-state energy (normalized) in a harmonic potential as a function of the statistics parameter $\alpha$ (fermions for $\alpha=0$), for different values of particle number $N$, and compared to the theoretical value $E^{\rm mTF}_{\alpha}/N^2$ in \eqref{eq:TF min}.}
\label{fig:num ener}
\end{center}
 \end{figure}

First, in Figure~\ref{fig:num ener} we plot the numerically evaluated energies for $N=25$, $N=50$ and $N=100$ particles in a harmonic trap, as a function of $\alpha$. We find a systematic numerical error compared to the theory~\eqref{eq:TF min}, but the expected oscillations are correctly reproduced, in particular the absolute maximum at $\alpha = 3/4$ expected from~\eqref{eq:TF constant} and Figure~\ref{fig:TF constant}.

Note that 
\begin{itemize}
  \item the agreement becomes noticeably better for larger $\alpha$, deeper in the mTF regime. This is also the most challenging situation, since the self-consistent magnetic field is larger. 
  \item the agreement improves when increasing $N$, as it should. This
    occurs rather slowly, and hence we have not pushed the numerics
    further than the already relatively costly $N=100$. 
\end{itemize}

\begin{figure}
\begin{center}
 \includegraphics[width=8cm]{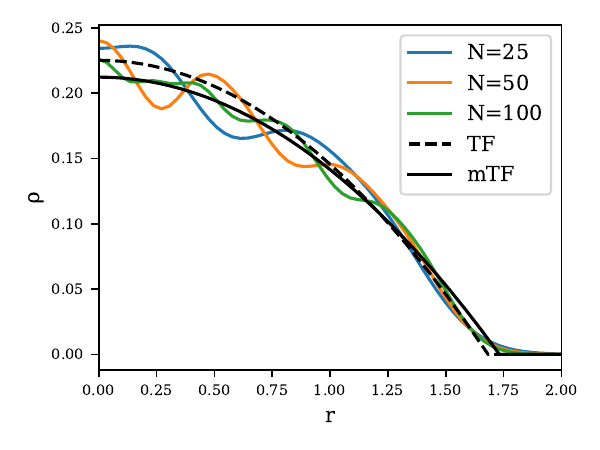}
\caption{Real/position space density $\rho(r)$ (normalized) in a harmonic trap for $\alpha=0.75$. Different values of $N$, and compared to the theoretical $\rho^{\rm TF}$ in \eqref{eq:TF dens} (dashed) resp. $\rho^{\rm mTF}_{\alpha}$ in \eqref{eq:TF min} (solid).}
\label{fig:density}
\end{center}
 \end{figure}

\begin{figure}
   \begin{center}
     \includegraphics[width=8cm]{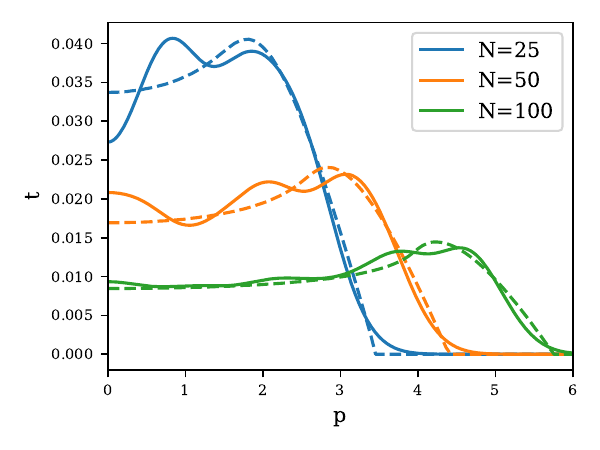} 
     \caption{Momentum space density $t(p)$ for $\alpha=3/4$, with $p$ rescaled by $N^{-1/2}$. Dashed lines
       are the semi-explicit theoretical guesses from \eqref{eq:mTF mom}--\eqref{eq:mTF A bis}.}
     \label{fig:t}
   \end{center}
 \end{figure}

Next, in Figures~\ref{fig:density} and \ref{fig:t}
we plot the position space and momentum space densities obtained for $N=25$, $N=50$ and $N=100$ particles in a harmonic trap~\eqref{eq:harm pot} at $\alpha = 3 / 4$. We again observe a correct match between numerical position space densities and mTF theory, but the oscillations of the former make it hard to argue for a clear signature of a $c(\alpha) \neq 1$. The trend with $N$ is however encouraging: around $N=100$, the mTF density (solid black line) starts to approximate the numerical plots better than the usual TF density (dashed black line).

The momentum space densities show, perhaps, our most interesting finding. There is a priori no obvious theory to compare them to, for mTF theory is not based on the usual ``position $\times$ momentum'' phase space. We however observe that some extrapolation of the TF momentum density fits the numerical curves pretty well. Namely, this is the semi-explicit formula~\eqref{eq:TF mom} where we replace $\rho^{\rm TF}$ by $\rho^{\rm mTF}_\alpha$ (given by~\eqref{eq:TF min}-\eqref{eq:TF constant} at $\alpha=3/4$). This leads to the expression
\begin{equation}\label{eq:mTF mom}
t^{\rm mTF}_\alpha (\bp)= \frac{1}{(2\pi)^2} \int_{\R^2}  \1 \left( \left|\bp + \alpha \bA_\alpha ^{\rm mTF} (\bx)\right|^2 \leq 4\pi \rho_\alpha^{\rm mTF} (\bx) \right) \rmd\bx
\end{equation}
with
\begin{equation}\label{eq:mTF A}
\bA_\alpha ^{\rm mTF} (\bx) 
:= \bA[\rho_\alpha^{\rm mTF}]
= \int_{\R^2} \frac{(\bx-\by)^\perp}{|\bx-\by|^2} \rho_\alpha ^{\rm mTF} (\by) \,\rmd\by.  
\end{equation}
and thus, if $V(\bx) = |\bx|^s$,
\begin{equation}\label{eq:mTF A bis}
\bA_\alpha ^{\rm mTF} (\bx) = \frac{N\bx^\perp}{4 c(\alpha)} \left( \left( 4c(\alpha) \frac{(s+2)}{s}\right)^{\frac{s}{s+2}} - \frac{2}{s+2} |\bx|^{s} \right) 
\end{equation}
on the support of $\rho^{\rm mTF}_\alpha$. Interestingly, we find that this use of $\rho^{\rm mTF}_\alpha$ instead of $\rho^{\rm TF}$ leads to better fits of the numerically, computed densities. These are shown in Figure~\ref{fig:t}. Note the non-monotonic shape of the curve, which marks the larger intensity of the effective magnetic field, and is a distinctive feature of the mTF regime.

 
\begin{figure}
\begin{center}
 \includegraphics[width=8cm]{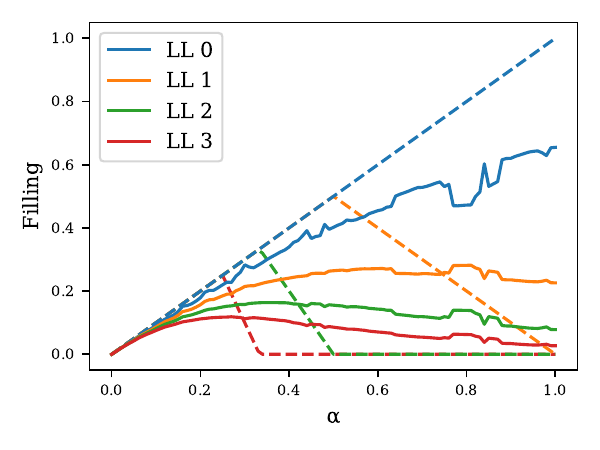} 
\caption{Relative filling $\mH (n,0)$ of local Landau levels at the origin, as a
  function of $\alpha$, normalized by $g_\varepsilon \ast \rho_{\gammaH} (0)$, compared to its
  theoretical value from \eqref{eq:local LL dist any} (dashed line). Parameters are $N=100$, $\epsilon = 0.1$,
  $V(\bx) = |\bx|^2$.}
\label{fig:LLs}
\end{center}
 \end{figure}

Finally, in Figure~\ref{fig:LLs}, we plot the local filling of Landau levels at the center of the trap, where the density is largest and hence mTF theory should be the most accurate. Namely, in accordance with~\eqref{eq:Husimi any}, we start from the numerically calculated density matrix $\gammaH$ to form the Husimi function
\begin{equation}\label{eq:num LL}
\mH (n,0) := \tr\left[ \Pi_{n,0} ^{B(0)} \, \gammaH \right],
\end{equation}
with $\Pi_{n,0} ^{B(0)}$ as in~\eqref{eq:coh op} (which seems to work better than~\eqref{eq:coh op 2} according to our numerical tests). To limit the influence of numerical errors we take the theoretical value
$$ B(0) = 2 \pi \alpha \rho^{\rm mTF}_{\alpha} (0).$$
There remains to choose the length $\varepsilon$ in~\eqref{eq:gaussian}. This is slightly tricky in practice, because at the values of $N$ we performed calculations for, there is not much room between magnetic length and system size to satisfy~\eqref{eq:lengths} (on top of that, the magnetic length seen as the radius of cyclotron orbits actually depends on the Landau level). A typical plot is shown in Figure~\ref{fig:LLs}, with comparison to the theoretical finding~\eqref{eq:local LL dist any}. 

Observe that 
$$ \mH (n,0) \leq \alpha \rho_{\gammaH} (0), \quad \sum_{n\geq 0} \mH (n,0) = g_\varepsilon \ast \rho_{\gammaH} (0), $$
so that it is fairly natural to normalize plots by the value 
$$\rho^{\rm mTF}_{\alpha} (0) \simeq g_\varepsilon \ast \rho_{\gammaH} (0)$$ 
and this is what we elected to do. Hence what is plotted is really the filling relative to the total density. 

We plot the (relative) filling of the lowest four Landau levels $n=0,1,2,3$. The fit with theory (as calculated from~\eqref{eq:local LL dist any}) is certainly not perfect, which can be related to discrepancies in energies seen on Figure~\ref{fig:num ener}. We do not observe the expected fact that any Landau level stays full until the next one is entirely empty. The noteworthy feature is however that starting around $\alpha\simeq 1/2$, a very limited number of Landau levels are required to account for the whole density. We interpret this as a transition towards mTF theory, whose main feature is indeed that only few Landau levels are taken into account. The residual Landau level mixing we observe may be due to finite size effects, variations of the magnetic length and the (chosen by hand) $\varepsilon$ length scale, etc. We have however observed a positive trend in $N$. For $N=100$ (shown in Figure~\ref{fig:LLs}), although the theoretical value in each level is not reproduced, we find that, for $\alpha \gtrsim 1/4$, at least $80\%$ of the total density is accounted for by the expected lowest $\lfloor \alpha^{-1} \rfloor + 1$ levels.
%

\section{Concluding Remarks}\label{sec:conclusions}

We have simulated a mean-field system of fermions coupled to magnetic flux tubes, effectively via a Chern-Simons self-generating gauge field, and compared the results to semi-classical approximations, based either on the usual ``position $\times$ momentum'' phase space or on the ``position $\times$ quantized cyclotron radius'' phase space. Being better suited for systems in large magnetic fields, we found that the latter approximation gives a better match to the numerics for general values of the magnetic flux, reproducing correctly ground-state energies and densities. A theoretical guess for the momentum density also fits numerics rather well, giving predictions on an observable which shows more marked dependence on the magnetic flux, and could thus serve as a better signature in putative experiments. Such a momentum density could indeed be obtained from time-of-flight measurements if the system under consideration is realized with cold atoms coupled to density dependent artificial gauge fields. This is one of the proposed realizations of anyons, quantum particles with fractional exchange statistics, already successfully implemented for 1D anyons~\cite{FroEtalTar-22,KwaEtalGre-24,DhaEtalNag-24}.

We note that our energy functional~\eqref{eq:Hartree} contains only the most salient feature one would expect from anyons or particles coupled to a Chern-Simons field, namely the self-consistent magnetic field~\eqref{eq:B Hartree}. It could be desirable to go beyond that by taking e.g. the following considerations into account:

\medskip

\noindent $\bullet$ If one directly inserts a Slater determinant trial state in the many-body energy functional~\eqref{eq:ener}, one finds as energy the functional~\eqref{eq:Hartree} that we studied thoroughly, but with the addition of extra terms. Besides exchange terms that would be of lower order for large $N$ one finds a two-body term~\cite[Section~2.1]{GirRou-21}
\begin{equation}\label{eq:twobody}
 \alpha^2 \iint_{\R^4} \frac{1}{|\bx - \by|^2} \left(\gamma (\bx,\bx)\gamma (\by,\by) - \left|\gamma (\bx,\by) \right|^2 \right) \rmd \bx \rmd \by. 
\end{equation}
Here $\gamma$ is, as in the main text, the one-particle density matrix of the Slater determinant, and the two terms inside the parenthesis are respectively the direct and exchange terms obtained from applying the Wick rule to compute the two-particle density matrix. A possible approximation of the above could be the ``direct only''
\begin{equation}\label{eq:twobody 2}
 \alpha^2 \iint_{\R^4} W (\bx - \by) \gamma (\bx,\bx)\gamma (\by,\by) \,\rmd \bx \rmd \by,
\end{equation}
where the effect of the exchange term (and thus the Pauli principle) is only taken into account by replacing $|\bx - \by|^{-2}$ with a $W(\bx-\by)$ regularized at short distances to make~\eqref{eq:twobody 2} finite.

We do not expect that the addition of terms of this type would alter our qualitative findings, while possibly blurring the main features somewhat, which is why we elected to stick with the pure Chern-Simons-Schr\"odinger formalism of Section~\ref{sec:Hartree}. 

\medskip

\noindent $\bullet$ Using more sophisticated trial states, incorporating repulsive Jastrow factors, would lead~\cite{AtaGirLun-25,ComCabOuv-91,Mashkevich-96,ComMasOuv-95,Ouvry-94} to the replacement of~\eqref{eq:twobody}-\eqref{eq:twobody 2} by a purely local term 
\begin{equation}\label{eq:twobody 3}
f (\alpha) \int_{\R^2} \gamma (\bx,\bx) ^2 \,\rmd \bx
\end{equation}
for some positive function $f(\alpha)$, see~\cite{AtaGirLun-25}. Such a term could have the intepretation of a Pauli term coupling the effective, density-dependent, magnetic field, to some spin of the particles. It could be incorporated with no extra difficulty in our basic scheme. With attractive Jastrow factors instead one could imagine a negative value of $f(\alpha)$.

\medskip

\noindent $\bullet$ Anyons realized as quasi-particles of condensed matter systems would have a non-zero extension, leading to a smearing of the magnetic flux over some, typically small, length scale. We have not considered this so far, but we do not expect this to influence our main conclusions dramatically, at least not at the level of a mean-field approximation. If anything, this would tend to smoothen the two-particle density of the system, making the mean-field approximation even more efficient. 

\medskip

\noindent $\bullet$ An additional, constant, external magnetic field could be added to the model, to better match the fractional quantum Hall effect context, and/or check for signs of superfluidity/superconductivity; see e.g.~\cite{HuMurRaoJai-21} and references therein.

\medskip

\noindent $\bullet$ As regards proposed implementations in cold atomic systems, the relationship between the gauge field and the matter density is in general more complicated than the pure flux-attachment~\eqref{eq:B Hartree}. We refer to the discussion in~\cite[Section~5]{ValWesOhb-20} for this aspect, that it could be desirable to take into account more precisely.

\bigskip

\noindent\textbf{Acknowledgments.} We thank Th\'eotime Girardot for discussions, the previous mathematical work on the almost-fermionic limit, and preliminary numerical work on the same topic~\cite{Girardot-21}.
D.L. also thanks Alireza Ataei and Dinh-Thi Nguyen for related discussions, 
as well as the Department of Mathematics of Politecnico di Milano for its kind hospitality during the spring of 2025 intensive period ``Quantum Mathematics at Polimi''.
Financial support from the Swedish Research Council 
(D.L., grant no.\ 2021-05328, ``Mathematics of anyons and intermediate quantum statistics'') 
is gratefully acknowledged.

\newpage

\end{document}